\documentstyle{amsart}
\input{bull-art}
\newtheorem{theorem}{Theorem}
\newtheorem{cor}{Corollary}

\newcommand{\bs}{\backslash}
\newcommand{\ds}{\displaystyle}
\newcommand{\fr}{\frac}
\newcommand{\hs}{\hspace}

\newcommand{\lt}{\left}

\newcommand{\ra}{\rightarrow}
\newcommand{\rt}{\right}

\newcommand{\beq}{\begin{equation}}
\newcommand{\eeq}{\end{equation}}
\newcommand{\Beq}{\begin{eqnarray*}}
\newcommand{\Eeq}{\end{eqnarray*}}
\newcommand{\BEQ}{\begin{eqnarray}}
\newcommand{\EEQ}{\end{eqnarray}}

\newcommand{\U}{\mbox{$\bold u$}}

\newcommand{\ld}{\lambda}
\newcommand{\Ome}{\Omega}

\newcommand{\R}{{\bf R}}

\newcommand{\sCoi}{C_{0}^{\infty}(\Ome)}

\newcommand{\sHo}{\widehat{H}_0^1 (\Ome)}
\newcommand{\sHoOm}{\widehat{H}_0^1 (\Ome_n )}

\newcommand{\sHoi}{H_0^1 (\Ome)}

\newcommand{\sLz}{L^2 (\Ome)}

\newcommand{\pa}{\partial}
\newcommand{\gr}{\nabla}
\newcommand{\la}{\Delta}

\newcommand{\Ugu}{\U\cdot\gr\U}
\newcommand {\gu}{\gr\U}

\newcommand{\sgu}{\gr u}
\newcommand {\Gu}{\,\|\gu}

\newcommand{\sGu}{\,\|\sgu}
\newcommand {\GA}{\|\gr\U_{0}\|}
\newcommand {\GU}{\Gu\|}
\newcommand{\sGU}{\sGu\|}

\newcommand {\lu}{\la\U}

\newcommand{\slu}{\la u}
\newcommand {\Lu}{\,\|\lu}

\newcommand{\sLu}{\,\|\slu}
\newcommand {\LU}{ \Lu\|}
\newcommand{\sLU}{\sLu\|}

\newcommand{\sUU}{\,\| u\|}

\newcommand{\ddt}{\fr{d}{dt}}
\newcommand{\dx}{\,dx}

\newcommand{\summ}{\sum_{n=1}^{m}}
\newcommand{\sumy}{\sum_{n=1}^{\infty}}
\newcommand{\supO}{\sup_{\Omega}}
\newcommand{\IO}{\int_{\Omega}}
\newcommand{\Iot}{\int_{0}^{t}}
\newcommand{\half}{\fr{1}{2}}

\begin{document}
\def\currentvolume{26}
\def\currentissue{2}
\def\currentyear{1992}
\def\currentmonth{April}
\def\copyrightyear{1992}
\def\currentpages{294-298}

\title[$L^2$-Laplacians on arbitrary domains]
{A sharp pointwise bound for functions 
\\ with $L^2$-Laplacians on arbitrary domains 
\\ and its applications}
\author{Wenzheng Xie}
\address{School of Mathematics, University of Minnesota, 
Minneapolis, 
Minnesota 55455}
\address{{\it E-mail address}:  {\defaultfont 
xie@@s5.math.umn.edu}}
\date{January 20, 1991 and, in revised form,
September 10, 1991}
\subjclass{Primary 26D10, 35B45, 35Q20}
\thanks{Partially supported by NSERC (Canada). The
contents of this paper have been presented to the Annual 
Meeting of the
American Mathematical Society, January 16--19, 1991}

\maketitle

\begin{abstract}
For all functions on an arbitrary open set 
$\Omega\subset\R^3$ 
with zero boundary values, we prove the optimal bound
\[
\sup_{\Omega}|u| \leq (2\pi)^{-1/2}
\left(\int_{\Omega}|\nabla u|^2 \,dx\,
\int_{\Omega}|\Delta u|^2 \,dx\right)^{1/4}.
\]
The method of proof is elementary and admits 
generalizations.
The inequality is applied to establish an existence 
theorem for the Burgers 
equation.
\end{abstract}

\section{Introduction}
In this note we announce the proof of the inequality
\begin{equation}
\sup_{\Omega}|u| \leq\frac{1}{\sqrt{2\pi}}
\left(\int_{\Omega}|\nabla u|^2 \,dx\,
\int_{\Omega}|\Delta u|^2 \,dx\right)^{1/4} 
\end{equation}
for functions with zero boundary values on 
three-dimensional domains.
The domain $\Ome$ can be any open set and the constant 
$1/\sqrt{2\pi}$ is
optimal. This best possible result is obtained by a new 
and elementary method, 
which is apparently also applicable to other elliptic 
operators. Thus, many 
known inequalities can be improved and new ones derived. 
Some of these will 
be given by the author in separate papers. Such 
inequalities  
are used in the study of nonlinear differential equations, 
see [1] and [2].

For smoothly bounded domains, 
one can combine the Sobolev inequality (see [3]) 
\[ \supO|u| \leq C_1 (\Ome)\|\gr u\|_{L^{2} (\Ome)}^{1/2} 
                               \|u\|_{H^{2} (\Ome)}^{1/2}\]
with the a priori estimate (see [4])
\begin{equation}
\|u\|_{H^2 (\Ome)} \leq C_2 (\Ome)\|\la u\|_{L^2 (\Ome)} 
\end{equation}
to obtain 
\begin{equation}
\supO|u| \leq C_3 (\Ome)\|\gr u\|_{L^{2} (\Ome)}^{1/2} 
                           \|\la u\|_{L^{2} (\Ome)}^{1/2} ,
\end{equation}
where $C_i (\Ome)$ are constants depending on the domain 
$\Ome$. 
However, the elliptic estimate (2)
fails to hold for domains with reentrant corners [5].

It was suggested to the author by Professor J. G. Heywood 
that (3) should be valid for nonsmooth domains as well, 
and its 
generalization to the Stokes operator would yield results 
for 
the Navier-Stokes equations in nonsmooth domains. 
Here, in \S4, we use (1) to derive a priori estimates 
and prove an existence theorem for the initial-boundary 
value problem of the 
Burgers equation with $H^1$ initial data, 
in an arbitrary open set, for the first time. 
\vspace{2mm}

\begin{center}{\sc 2. The main results} 
\end{center}\vspace{2mm}

Let $\,\Ome\,$ be an arbitrary open set in $\,\R^3 \,$. 
Let $\,\|\cdot\| \,$ denote the $\,\sLz \,$ norm.
The homogeneous Sobolev space $\sHo\,$ is defined to be 
the completion of $\,\sCoi \,$ in the Dirichlet norm 
$\,\|\gr\cdot\| \,$, where $\,\gr\,$ is the gradient. 
Let $\,\la\,$ denote the Laplacian in the sense of 
distributions. 
Our main result is

\begin{theorem} 
For all $\, u\in\sHo \,$ with $\,\slu\in\sLz \,$, there 
holds 
\[\supO|u| \leq\fr{1}{\sqrt{2\pi}}\sGU^{1/2}\sLU^{1/2} .\]
The constant $\,1/\sqrt{2\pi} \,$ is optimal for each 
$\Ome$.
\end{theorem} 

The space $\sHo\,$ contains the standard Sobolev space 
$\,\sHoi\,$. It 
contains functions that are not square integrable for some 
unbounded domains.
If $\,u\in\sHoi\,$, then using $\sGU^2 =-\IO u\la u \dx 
\leq\sUU\sLU$, we 
obtain 

\begin{cor}
If $\,u\in\sHoi \,$ and $\,\slu\in\sLz \,$, then $\,u\,$ 
also satisfies 
\[ \supO|u|\leq\fr{1}{\sqrt{2\pi}}\sUU^{1/4}\sLU^{3/4} .\]
\end{cor}

In particular, we obtain a pointwise bound for any 
normalized
eigenfunction of the Laplacian, in terms of its 
corresponding eigenvalue.

\begin{cor}
If $\,u\,$ satisfies 
\[ -\la u =\ld u\,, \hs{6mm}u\in\sHoi\,, 
\hs{6mm}\|u\|=1\,,\] 
then 
\[ \supO|u|\leq\fr{\ld^{3/4}}{\sqrt{2\pi}} \,. \]
\end{cor}

All of the above results are also valid for vector-valued 
or complex-valued 
functions. The constants in the corollaries, however, are 
not optimal.

\vspace{2mm} \begin{center}{\sc 3. Outline of proof} 
\end{center}

The proof of (1) has four steps. 

{\em Step} 1. 
First we assume that $\,\Ome\,$ is bounded, with a 
$\,C^{\infty} \,$ boundary 
$\,\pa\Ome\,$. It is well known that there exist 
eigenfunctions 
$\,\{\phi_{n}\}\,$ of the Laplacian
that form a complete orthonormal basis of $\,\sLz \,$,
satisfying
\[ -\la \phi_{n}=\ld_{n}\phi_{n} \,, 
\qquad \phi_{n} |_{\pa\Ome}=0 \,,\]
where $\,\ld_n >0\,$ are the eigenvalues, $n=1,2,\ldots\,$. 

Let $\,x_0 \in\Ome\,$ and $\,m\geq 1\,$ be fixed. 
For functions of the form $u(x)=\summ  c_{n} 
\phi_{n}(x)\,$, we have
\[ \fr{u^2 (x_0 )}{\sGU\sLU} 
=\fr{\ds \lt(\sum\nolimits^m_{n=1}      c_{n}\phi_{n}(x_0 
)\rt)^2 }
    {\ds \left(\sum\nolimits^m_{n=1} \ld_{n}    
c_{n}^{2}\right)^{1/2} 
         \left(\sum\nolimits^m_{n=1} 
\ld_{n}^{2}c_{n}^{2}\right)^{1/2} } \,.\]
This quotient is a smooth and homogeneous function of 
$\,(c_1 ,\dotsc,c_m )\,$ 
in $\,\bold R^m\backslash\{0\}$. Hence, at some point 
$\,(\tilde{c}_1 ,\dotsc, \tilde{c}_m ) \,$, 
it attains its maximum value.  The maximum value can be 
written as
\[ 4\sqrt{\mu}\summ\left( \fr{\phi_n (x_0 )}{\mu+
\ld_n}\right)^2 \,,\]
where $\mu=\summ \ld_n^2 \tilde{c}_n^2 /\summ \ld_n 
\tilde{c}_n^2 $.

{\em Step} 2. We introduce the Green function 
$G(x;x_0 ,\mu)$ for the Helmholtz equation
\[ \la G =\mu G -\delta(x-x_0 )\,,
\hspace{7mm}  G|_{\pa\Ome} =0 \,.\] 
By the maximum principle, we have
\[ 0\leq G(x;x_0 ,\mu)\leq \fr{e^{-\sqrt{\mu}|x-x_0 
|}}{4\pi |x-x_0 |},\quad 
\forall x\in\Ome\bs\{x_0 \}\,,\]
the upper bound being the fundamental solution. Hence
\[ \IO G^2\dx \leq  
\int_{0}^{\infty}\left(\fr{e^{-\sqrt{\mu}r}}{4\pi 
r}\right) ^2 4\pi r^2 \,dr 
=\fr{1}{8\pi\sqrt{\mu}} .\]
By Parseval's equality and Green's formula, we have
\[ \IO G^2\dx=\sumy\lt(\IO G\phi_n\dx\rt)^2
=\sumy\lt(\fr{\phi_n (x_0 )}{\mu+\ld_n}\rt)^2  .\]
Therefore
\[ \fr{u^2 (x_0 )}{\sGU\sLU} 
\leq 4\sqrt{\mu}\summ\left( \fr{\phi_n (x_0 )}{\mu+
\ld_n}\right)^2 
\leq 4\sqrt{\mu}\IO G^2\dx 
\leq\fr{1}{2\pi} \,.\]
Thus, (1) is true for any 
function of the form $u(x)=\summ  c_{n} \phi_{n}(x)\,$.

{\em Step} 3.  Now, let $\,u\,$ be any function in $\,\sHo 
\,$ such that 
$\,\slu\in\sLz\,$. 
Let $\,u_n \,$ be the projection of $\,u\,$ in 
span$\{\phi_1 ,\dotsc,\phi_n \}$. We have
$\,\sGu_n \| \leq\sGU \,$, $\,\sLu_n \| \leq\sLU\,$, and   
$\,\lim_{n\ra\infty}u_n =u$ in $\,\sLz\,$. 
It follows that (1) remains valid.

{\em Step} 4.  Now we proceed to prove Theorem 1.  
We can choose a sequence of bounded domains $\,\Ome_n \,$ 
with smooth 
boundaries such that 
$\ds\Ome_1 \subset\Ome_2 \subset\cdots$ and 
$\,\bigcup_{n=1}^{\infty}\Ome_n =\Ome \,$.
For each $\,n\geq 1 \,$, 
there exists a unique $\,u_n \in\sHoOm \,$ such that 
\[ \int_{\Ome_n}\sgu_n \cdot\gr 
v\dx=\int_{\Ome_n}\sgu\cdot\gr v\dx,\quad
\forall v\in\sHoOm  \,,\]
by the Riesz representation theorem. From this we obtain
$\,\sGu_n \|_{L_{2} (\Ome_n )} \leq \sGU \,$, $\,\slu_n 
=\slu|_{\Ome_n}\,$, 
and $\,\lim_{n\ra\infty}u_n =u\,$ in $\,\sHo\,$, hence in 
$\,L^6 (\Ome)\,$. 
Therefore the proof of (1) is completed.

Let $u(x)=(1-e^{-|x|})/|x|$; then we have 
\[ \sup_{x\in\bold R^3 }|u(x)|= 1\,,\hs{6mm}
\int_{\bold R^3}|\gr u|^2\, \dx = 2\pi\,,\hs{6mm} 
\int_{\bold R^3}|\la u|^2\, \dx = 2\pi\,.\]
Hence the equality in (1) holds for $u$.  
By cutting-off $u$, we explicitly construct a sequence of 
functions 
$u_n$ with compact support such that
\[ u_n (0)\ra 1\,,\hs{6mm}
\int_{\bold R^3}|\gr u_n |^2\, \dx \ra 2\pi\,,\hs{6mm} 
\int_{\bold R^3}|\la u_n |^2 \,\dx \ra 2\pi\,,\]
as $\,n\ra\infty\,$. 
Given any open set $\,\Ome\,$, by scaling, we obtain a new 
sequence of 
functions with compact support in $\Ome$. Since the product 
\[
\|\gr u_n \|_{L^2 (\bold R^3 )}^{1/2}\,\|\la u_n \|_{L^2 
(\bold R^3 )}^{1/2}\]
is scale invariant, it is seen that the constant 
$1/\sqrt{2\pi}$ in (1) 
is the best possible.

\begin{center}{\sc 4. Application to the Burgers equation} 
\end{center}

The time-dependent Burgers equation 
\begin{equation}
\fr{\pa \U}{\pa t} +\Ugu =\nu\la\U 
\end{equation}
is sometimes studied for its analogy with the 
Navier-Stokes equations, 
with the three-dimensional vector-valued function 
$\U\!(x,t)$ 
representing the velocity field and the 
positive constant $\nu$ the viscosity coefficient. 
We consider as spatial domain an arbitrary open set in 
$\bold R^3$
and seek $\U$ that vanishes on the boundary and takes 
an initial value $\U_0 \in\sHo^3$.
From (4) and the vector version of (1), we have
\Beq
\half\ddt\GU^2+ \nu\LU^2 &=& \IO \Ugu\cdot\lu\dx \\
&\leq& \sup |\U|\GU\LU   \\
&\leq& \ds\fr{1}{\sqrt{2\pi}}\GU^{3/2}\LU^{3/2} \,.
\Eeq
By using Young's inequality and a comparison theorem, 
we obtain
\[ \Gu(t)\|^2 \leq\fr{\GA^2}{\sqrt{1-t/T} } \,,\]
and 
\[ \Iot\Lu(s)\|^2 \,ds \leq 
\fr{\GA^2}{2\nu\lt(1-\sqrt[6]{t/T}\rt)\sqrt{1-\sqrt{t/T}}} 
\,, \]
for $0\leq t<T$, where 
\[ T=\fr{256\pi^2 \nu^3 }{27\GA^4} \,.\] 
Beginning with these a priori estimates, and using the 
methods of [1] and 
[6], the following theorem is established [7].

\begin{theorem} 
For any open $\Ome\subset\bold R^3$ and any $\U_0 
\in\sHo^3$, 
there exists a unique function
\[ \U\in C([0,\,T),\,\sHo )^3 \cap\, 
C^{\infty}(\Ome\times(0,\,T))^3 
\,\cap\, C^{\infty}((0,\,T),\,L_{\infty}(\Ome))^3 \,,\] 
satisfying the Burgers equation {\em (4)} and taking the 
initial value $\U_0$. 
\end{theorem} 

\begin{center}\sc Acknowledgments \end{center}

I am grateful to J. G. Heywood for suggesting the topic 
for this paper and for his helpful advice. 
I would also like to thank L. Rosen for helpful 
discussions. 



\end{document}